\documentstyle{article}
\begin{document}

\newcommand{\HDS}{\vrule width0pt height2.3ex depth1.05ex\displaystyle}

\def\f#1#2{{{\HDS #1}\over{\HDS #2}}}

\def\lpravilo#1{ \makebox[-.5em][r]{\mbox{\it #1}} {\mbox{\hspace{0.5em}}}}

\def\str{\rightarrow}

\newcommand{\qed}{\hfill $\Box$}

\def\k#1#2{\stackrel{\raisebox{-2pt}{\mbox{\tiny $#1$}}}{k}^
{\raisebox{-8pt}{\scriptsize #2}}}

\def\G{\mbox{\it Gen}}

\def\R{\mbox{\it Rel}}

\def\mj{\mbox{$\boldmath 1$}}

\def\cl{\mbox{\rm Cl}}

\title{\bf Generality of Proofs \\ and its Brauerian Representation}

\author{{\sc Kosta Do\v sen} and {\sc Zoran Petri\' c}
\\[.05cm]
\\Mathematical Institute, SANU \\
Knez Mihailova 35, p.f. 367 \\
11001 Belgrade, Yugoslavia \\
email: \{kosta, zpetric\}@mi.sanu.ac.yu}
\date{}
\maketitle

\begin{abstract}
\noindent The generality of a derivation is an equivalence relation on
the set of occurrences of variables in its premises and conclusion such
that two occurrences of the same variable are in this relation if and
only if they must remain occurrences of the same variable in every
generalization of the derivation. The variables in question are
propositional or of another type. A generalization of the derivation
consists in diversifying variables without changing the rules of
inference.

This paper examines in the setting
of categorial proof theory the conjecture that two derivations
with the same premises and conclusions stand for the
same proof if and only if they have the same generality. For that purpose
generality is defined within a category whose arrows are
equivalence relations on finite
ordinals, where composition is rather complicated.
Several examples are given of deductive systems of derivations covering
fragments of logic, with the associated map into the category of
equivalence relations of generality.

This category is isomorphically represented in the category whose arrows are
binary
relations between finite ordinals, where composition is the usual
simple composition of relations.
This representation is related to a classical representation result of
Richard Brauer.

\vspace{0.3cm}

\end{abstract}

\section{Introduction}

Up to now, two conjectures
were made concerning identity criteria for proofs.
The first, and dominating conjecture, which we call the {\it Normalization
Conjecture}, was made by Prawitz in \cite{P71}. There Prawitz raised the
question of identity criteria for proof as one of the central questions
of general proof theory, and conjectured that two derivations in
natural deduction should stand for the same proof iff they both reduce to
the same normal form. (For details and further references concerning
the Normalization Conjecture see \cite{D02}.) Via the Curry-Howard
correspondence, the Normalization Conjecture says that identity of
proofs is caught by equality between typed lambda terms, and this shows
that for fragments of intuitionistic logic the equivalence relation
on derivations envisaged by Prawitz is a well-grounded mathematical
notion. A further corroboration for this conjecture is brought by
the categorial equivalence between typed lambda calculuses and
cartesian closed categories, which was demonstrated by Lambek
(see \cite{L74} and \cite{LS86}, Part I). Cartesian closed categories
provide an alternative description of proofs in intuitionistic
conjunctive-implicative logic, and the equivalence relation between
proofs given by identity of arrows in cartesian closed categories is in
accordance with the Normalization Conjecture.

The second conjecture concerning identity criteria for proofs, which we call
the {\it Generality Conjecture}, was suggested by Lambek in
\cite{L68} (pp. 287, 316) and \cite{L69} (p. 89). According to
this conjecture, two derivations with the
same premises and conclusions should stand for the same proof iff for
every generalization of one of them, which consists in diversifying
the variables without changing the rules of inference, there is a
generalization of the other derivation so that in the two generalizations
we have the same premises and conclusions. For example, let
$\k{\wedge}{1}_{p,p}:p\wedge p\str p$ and
$\k{\wedge}{2}_{p,p}:p\wedge p\str p$
be obtained respectively by the first rule and second rule of conjunction
elimination, let $\langle \k{\wedge}{1}_{p,p},
\k{\wedge}{2}_{p,p}\rangle:p\wedge p\str
p\wedge p$ be obtained from $\k{\wedge}{1}_{p,p}$ and
$\k{\wedge}{2}_{p,p}$ by
conjunction introduction, and let $1_{p\wedge p}:p\wedge p\str p\wedge p$
be the identity derivation from $p\wedge p$ to $p\wedge p$. For every
generalization of $\langle \k{\wedge}{1}_{p,p},
\k{\wedge}{2}_{p,p}\rangle$---and there
is only one sort of generalization here, namely,
$\langle\k{\wedge}{1}_{q,r},
\k{\wedge}{2}_{q,r}\rangle:q\wedge r\str q\wedge r$---we have
a generalization of $1_{p\wedge p}$, namely,
$1_{q\wedge r}:q\wedge r\str q\wedge r$.
So, by the Generality Conjecture,
$\langle \k{\wedge}{1}_{p,p},\k{\wedge}{2}_{p,p}\rangle$ and
$1_{p\wedge p}$ should stand
for the same proof, which is in accordance with the Normalization Conjecture
in its beta-eta version. On the other hand, $\k{\wedge}{1}_{p,p}$ can be
generalized to $\k{\wedge}{1}_{q,r}:q\wedge r\str q$,
and there is no generalization
of $\k{\wedge}{2}_{p,p}$ of the type $q\wedge r\str q$. So, by the Generality
Conjecture, $\k{\wedge}{1}_{p,p}$ and $\k{\wedge}{2}_{p,p}$
would not stand for the same proof;
this is again in accordance with the Normalization Conjecture.
(We will see in Section 3 that when we go beyond the conjunctive fragment
of logic, the two conjectures need not accord with each other any more.)

The notion of generality of derivation, which is implicit in the Generality
Conjecture, will be made precise in the context of category theory. A proof
will be an arrow in a category, and a derivation will be represented by an
arrow term. The source of the arrow is the premise and the
target the conclusion of the proof, and analogously for arrow terms and
derivations. The generality of an arrow term will be an equivalence relation
on the set of occurrences of variables in the source and target of this
arrow term. If two occurrences in this set are equivalent, then they
must be occurrences of the same variable (but not vice versa). The
Generality Conjecture will say that two arrow terms with the same
source and target stand for the same arrow iff they have the same
generality.

Note that Lambek's formal notion of generality
from \cite{L68} and \cite{L69}, with which he became dissatisfied in
\cite{L72} (p. 65), is not the notion we will define
in this paper. We take from Lambek just the inspiration of his
intuitive remarks about generality of proofs.

In the next section, we will introduce a category called \G, in
which we will define formally the generality of derivations.
We may conceive this category as obtained by abstraction from the category
whose arrows are proofs and whose objects are formulae by replacing formulae
with sequences of
occurrences of variables. From arrows that stood for proofs we keep
only the equivalence relation holding between occurrences of
variables that
must in every generalization remain occurrences of the same variable.
Since sequences of occurrences of variables are finite we represent
them by finite ordinals. The variables that motivate the introduction
of \G\ are of an arbitrary type: they may be propositional
variables, or individual variables, or of any other type of
variable in a language from which we take our formulae. We
may also deal simultaneously with several types of variable.

In Section 3, we will give several examples of deductive systems,
covering fragments of logic, and of generality obtained
via maps into \G\ defined for these deductive systems, which induces identity
of proofs according to the Generality Conjecture.

The definition of composition of arrows in \G\ is rather involved, and we
will not attempt to give a direct formal proof of the associativity of this
composition. An indirect proof of this associativity will be provided
by an isomorphic representation of \G\ in the category \R\ whose arrows are
binary
relations between finite ordinals, where composition is just
composition of relations. We devote Sections 4-6 to this matter. In Section
4 we first deal with some very general questions about representing
equivalence relations, while in Section 5 we give our representation
of \G\ in \R, and prove the isomorphism of this representation.
In Section 6 we show that our representation of \G\ in \R\ is
related to a classical result of representation theory
due to Richard Brauer.

The representation of \G\ in \R\ does not serve only to prove the
associativity of the composition of \G, but it may be of an independent
interest. Although in passing from \G\ to \R\ we have a functor
involving an exponential function, \R\ may be helpful in
computations, since composing binary relations is like multiplying
matrices. The representation of \G\ in \R\ sheds also some light on
Brauer's representation result, which it generalizes to a certain
extent.

In a companion to this paper \cite{DP03} we will consider
a category more general than \G\ obtained by replacing equivalence relations
by preordering, i.e. reflexive and transitive, relations.
This category will also be isomorphically represented in \R,
in a manner that generalizes further Brauer's representation.

\section{The category \G}

The objects of the category \G\ are finite ordinals, and its arrows
$R:n\str m$ are equivalence relations $R\subseteq (n+m)^2$.
As usual, we write sometimes $xRy$ for $(x,y) \in R$.

If $X\subseteq \omega$, then $X^{+k}=_{\mbox{\scriptsize\it def}}
\{n+k\mid n\in X\}$. For $n,m\in \omega$,
the sum $n+m$ is by definition $n\cup m^{+n}$. For an arbitrary
relation $R\subseteq (n\cup m^{+n})^2$ and $k\in\omega$, let $R^{+k}
\subseteq (n^{+k}\cup m^{+n+k})^2$ be the relation defined by
\[x R^{+k} y \quad {\mbox{\rm iff}}\quad (x-k) R (y-k).\]
Note that if $R$ is an equivalence relation in \G, then for $k\neq 0$ the
equivalence relation $R^{+k}$ is not an arrow of \G.

The identity arrow $\mj_n$ in \G\ is the equivalence relation
$\mj_n\subseteq (n+n)^2$ defined by

\[(x,y) \in \mj_n  \quad {\mbox{\rm iff}}
\quad x=y\;({\makebox{\rm mod }}n).\]

We define composition of arrows in \G\ in the following manner.
Let $\cl(R)$ be the transitive closure of the binary relation $R$.
Suppose
we have the equivalence relations $R_1:n\str m$ and $R_2:m\str k$; that is,
$R_1\subseteq(n\cup m^{+n})^2$ and $R_2\subseteq(m\cup k^{+m})^2$.
Then
$R_1\cup R_2^{+n}\subseteq (n\cup m^{+n}\cup k^{+m+n})^2$,
and it is easy to check that
$\cl(R_1\cup R_2^{+n})$ is an equivalence relation.

For $x\in n\cup k^{+n}$ and $m\in\omega$ let
\[x'=\left\{\begin{array}{lll}
x & {\mbox{\rm if}} & x\in n
\\
x+m & {\mbox{\rm if}} & x\in k^{+n}.
\end{array}
\right.
\]
Then for an arbitrary relation $R\subseteq (n\cup k^{+m+n})^2$ let
$R^\dagger\subseteq(n\cup k^{+n})^2$ be defined by
\[ x R^\dagger y \quad {\mbox{\rm iff}} \quad x' R y'.\]
It is easy to see that if $R$ is an equivalence relation, then
so is $R^\dagger$.

Then we define the composition $R_2\ast R_1:n\str k$ of the arrows
$R_1$ and $R_2$ above by
\[
R_2\ast R_1 =_{\mbox{\scriptsize\it def}} ((\cl(R_1\cup R_2^{+n})\cap
(n\cup k^{+m+n})^2)^\dagger.
\]

It can be checked that $R_2\ast R_1$ is indeed an equivalence relation on
$n\cup k^{+n}$. It is more involved to check formally that $\ast$ is
associative, though this is rather clear if we represent this
operation geometrically (as this is done in \cite{EK66}, and also
in categories of tangles; see \cite{K95}, Chapter 12, and references therein).
For example, for $R_1\subseteq (3\cup 9^{+3})^2=12^2$ that
corresponds to the partition
$\{\{0,3\}$, $\{4,5\}$, $\{1,6\}$, $\{7,8\}$, $\{2,9\}$, $\{10,11\}\}$, and
$R_2\subseteq (9\cup 1^{+9})^2=10^2$ that corresponds to the partition
$\{\{0,1\}$, $\{2,9\}$, $\{3,4\}$, $\{5,6\}$, $\{7,8\}\}$, the composition
$R_2\ast R_1\subseteq (3\cup 1^{+3})^2=4^2$ that corresponds
to the partition $\{\{0,3\}$, $\{1,2\}\}$ is obtained from the following
drawing

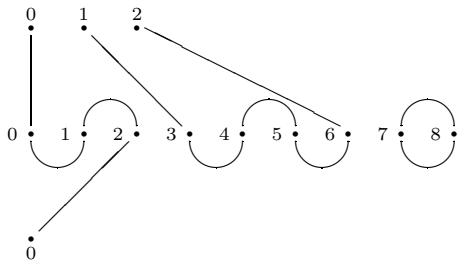
\begin{center}
\begin{picture}(160,120)

\put(3,23){\line(1,1){34}}
\put(0,63){\line(0,1){34}}
\put(57,63){\line(-1,1){34}}
\put(117,63){\line(-2,1){74}}

\put(0,20){\circle*{2}}
\put(0,60){\circle*{2}}
\put(20,60){\circle*{2}}
\put(40,60){\circle*{2}}
\put(60,60){\circle*{2}}
\put(80,60){\circle*{2}}
\put(100,60){\circle*{2}}
\put(120,60){\circle*{2}}
\put(140,60){\circle*{2}}
\put(160,60){\circle*{2}}
\put(0,100){\circle*{2}}
\put(20,100){\circle*{2}}
\put(40,100){\circle*{2}}

\put(10,57){\oval(20,20)[b]}
\put(30,63){\oval(20,20)[t]}
\put(70,57){\oval(20,20)[b]}
\put(110,57){\oval(20,20)[b]}
\put(90,63){\oval(20,20)[t]}
\put(150,57){\oval(20,20)[b]}
\put(150,63){\oval(20,20)[t]}

\put(0,17){\makebox(0,0)[t]{\scriptsize$0$}}
\put(-5,60){\makebox(0,0)[r]{\scriptsize$0$}}
\put(15,60){\makebox(0,0)[r]{\scriptsize$1$}}
\put(35,60){\makebox(0,0)[r]{\scriptsize$2$}}
\put(55,60){\makebox(0,0)[r]{\scriptsize$3$}}
\put(75,60){\makebox(0,0)[r]{\scriptsize$4$}}
\put(95,60){\makebox(0,0)[r]{\scriptsize$5$}}
\put(115,60){\makebox(0,0)[r]{\scriptsize$6$}}
\put(135,60){\makebox(0,0)[r]{\scriptsize$7$}}
\put(155,60){\makebox(0,0)[r]{\scriptsize$8$}}
\put(0,103){\makebox(0,0)[b]{\scriptsize$0$}}
\put(20,103){\makebox(0,0)[b]{\scriptsize$1$}}
\put(40,103){\makebox(0,0)[b]{\scriptsize$2$}}

\end{picture}
\end{center}

\noindent In this example all members of partitions have two elements, but
this is by no means necessary.

The complexity of the definition of $\ast$ in \G\ motivates the introduction
of an isomorphic representation of \G, considered in Section 5, in
which composition will be defined in an elementary way. The proof that
we have
there an isomorphic representation will provide an indirect proof
that $\ast$ in \G\ is associative, and that \G\ is indeed a category.

\section{Generality of derivations in fragments of logic}

The language $\cal L$ of conjunctive logic is built from a
nonempty set of propositional variables $\cal P$ with the binary connective
$\wedge$ and the propositional constant, i.e. nullary connective, $\top$
(the exact cardinality of $\cal P$ is not important here). We use the
schematic letters $A,B,C,\ldots$ for formulae of $\cal L$.

We have the following axiomatic derivations for every $A$ and $B$ in
$\cal L$:
\[
\begin{array}{l}
1_A:A\str A,
\\
\k{\wedge}{1}_{A,B}:A\wedge B\str A,
\\
\k{\wedge}{2}_{A,B}:A\wedge B\str B,
\\
\k{\wedge}{}_A:A\str\top,
\end{array}
\]
and the following inference rules for generating derivations:
\[
\f{f:A\str B \quad g:B\str C}{g\circ f:A\str C}
\]
\[
\f{f:C\str A \quad g:C\str B}{\langle f,g\rangle:C\str A\wedge B}
\]
This defines the deductive system $\cal D$ of conjunctive logic (both
intuitionistic and classical). In this system $\top$ is included as an
``empty conjunction''. For $f:A\str B$ and $g:C\str D$ let $f\wedge
g:A\wedge C\str B\wedge D$ abbreviate $\langle
f\circ\k{\wedge}{1}_{A,C}, g\circ\k{\wedge}{2}_{A,C}\rangle$.

We define now a function $G$ from $\cal L$ to the objects of \G\ by
taking that $G(A)$ is the number of occurrences of propositional
variables in $A$. Next we define inductively a function, also denoted
by $G$, from the derivations of $\cal D$ to the arrows of \G:
\[
\begin{array}{ll}
G(1_A)=\mj_{G(A)}, &
\\[.1cm]
(x,y)\in G(\k{\wedge}{1}_{A,B})  & {\mbox{\rm iff}} \quad
x=y\; ({\mbox{\rm mod }}G(A\wedge B)),
\\[.1cm]
(x,y)\in G(\k{\wedge}{2}_{A,B}) & {\mbox{\rm iff}} \quad
(x,y < G(A) \;{\mbox {\it and }} x=y)
\;{\mbox {\it or }}
(x,y \geq G(A) \;{\mbox {\it and }}
\\ & \hspace{5cm} x=y\; ({\mbox{\rm mod }} G(B))),
\\[.1cm]
(x,y)\in G(\k{\wedge}{}_{A}) & {\mbox{\rm iff}} \quad x=y,
\\[.2cm]
{\makebox[1em][l]{$G(g\circ f)=G(g)\ast G(f);$}} &
\end{array}
\]
for $x\in G(C)\cup G(B)^{+G(A)+G(C)}$ let
\[
x''=\left\{\begin{array}{ll}
x & {\mbox{\rm if }} x\in G(C)
\\
x-G(A) & {\mbox{\rm if }} x\in G(B)^{+G(A)+G(C)},
\end{array}
\right.
\]
and let $(x,y)\in G(g)^{\dagger\dagger}$ iff $(x'',y'')\in G(g)$; then
\[
\makebox[30em][l]{$G(\langle f,g\rangle)=\cl (G(f)\cup G(g)^
{\dagger\dagger}).$}
\]

Then, in accordance with the Generality Conjecture, we stipulate that
two derivations $f,g:A\str B$ of $\cal D$ are equivalent iff
$G(f)=G(g)$. A proof of $\cal D$ is then the equivalence class of a
derivation with respect to this equivalence relation.

It is then clear that proofs make the arrows of a category $\cal C$, with
the obvious sources and targets,
and that $G$
gives rise to a faithful functor from $\cal C$ to \G. The category
$\cal C$ happens to be the free cartesian category generated by the set of
propositional variables $\cal P$ as the generating set of objects (this
set may be conceived as a discrete category).
This fact about $\cal C$ follows from the coherence result for cartesian
categories treated in \cite{DP01} and \cite{P02}.
(Cartesian categories are
categories with all finite products, including the empty product, i.e.
terminal object. The category $\cal C$ can be equationally presented;
see \cite{LS86}, Chapter I.3, or \cite{DP01}.)

The language of disjunctive logic is dual to the language we had above:
instead of $\wedge$ and $\top$
we have $\vee$ and $\bot$ in $\cal L$, and instead
of $\k{\wedge}{\it i}$,
$\k{\wedge}{}$ and $\langle\; ,\; \rangle$ we have
in the corresponding deductive system $\cal D$
\[
\begin{array}{l}
\k{\vee}{1}_{A,B}:A\str A\vee B,
\\
\k{\vee}{2}_{A,B}:B\str A\vee B,
\\
\k{\vee}{}_A:\bot\str A,
\end{array}
\]
\[
\f{f:A\str C \quad g:B\str C}{[ f,g]:A\vee B\str C}
\]
Dually to what we had above, we obtain that equivalence of derivations
induced by generality makes of $\cal D$ the free category with all
finite coproducts generated by $\cal P$.

Let us now assume we have in $\cal D$ both $\wedge$ and $\vee$, but
without $\top$ and $\bot$, and let us introduce the category of proofs
$\cal C$ induced by generality as we did above. The category $\cal C$
will, however, not be the free category with nonempty finite products and
coproducts generated by $\cal P$, in spite of the coherence result of
\cite{DP02}. The problem is that in $\cal C$ we don't have the equation
\[\k{\wedge}{1}_{A,B}\circ(1_A\wedge[1_B,1_B])=\;\k{\wedge}{1}_{A,B\vee B}.\]
The graphs with respect to which coherence is proved in \cite{DP02}
don't correspond to the equivalence relations of \G. The graphs of
\cite{DP01}, with respect to which coherence can be proved for
cartesian categories are also not the equivalence relations of \G, but
in this case one easily passes to \G\ from the category whose arrows
are these graphs.

The coherence result of \cite{KML71}
for symmetric monoidal closed categories without
the unit object I can be used to show that generality in
the appropriate deductive system,
which corresponds to the
implication-tensor fragment of intuitionistic linear logic,
gives rise to the free symmetric
monoidal closed category without I. The graphs
of \cite{KML71} correspond to the equivalence relations of \G.

If we add intuitionistic implication to conjunctive or
conjunctive-disjunctive
logic, the equivalence of derivations induced by generality does not
give rise to free cartesian closed or free bicartesian
closed categories (for counterexamples see \cite{DP02}, Section 1, and
\cite{S75}).

Up to now, we considered generality of derivations in fragments of
logic with respect to propositional variables. In the next example this
generality will be considered with respect to individual variables.

Let now $\cal L$ be the language generated from a nonempty set of
individual variables $\{x,y,z,\ldots\}$, whose exact
cardinality is again not important, with the binary relational symbol
$=$, the binary connective $\wedge$ and the propositional constant $\top$.

As axioms and inference rules for our deductive system $\cal D$ we now
have whatever we had for conjunctive logic in the preceding section
plus all axioms of the form
\vspace{.3cm}
\[
\begin{array}{l}
r_x:\top\str x=x,
\\
s_{x,y}:x=y\str y=x,
\\
t_{x,y,z}:x=y\wedge y=z\str x=z.
\end{array}
\]
\vspace{.2cm}

The function $G$ from $\cal L$ to the objects of \G\ is now defined  by
taking that $G(A)$ is the number of occurrences of individual variables
in $A$. We define inductively the function $G$ from the derivations of
$\cal D$ to the arrows of \G\ as we did in the preceding section for
conjunctive logic with the following additional clauses, where $X^+$ is
the reflexive and symmetric closure of $X$,
\[
\begin{array}{l}
G(r_x)=\{(0,1)\}^+,
\\
G(s_{x,y})=\{(0,3),(1,2)\}^+,
\\
G(t_{x,y,z})=\{(0,4),(1,2),(3,5)\}^+.
\end{array}
\]
Graphically, these clauses correspond to

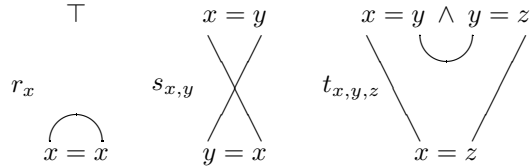
\begin{center}
\begin{picture}(180,80)

\put(60,20){\line(1,2){20}}
\put(80,20){\line(-1,2){20}}
\put(140,20){\line(-1,2){20}}
\put(160,20){\line(1,2){20}}

\put(10,20){\oval(20,20)[t]}
\put(150,60){\oval(20,20)[b]}

\put(10,17){\makebox(0,0)[t]{$x=x$}}
\put(70,17){\makebox(0,0)[t]{$y=x$}}
\put(150,17){\makebox(0,0)[t]{$x=z$}}
\put(70,65){\makebox(0,0)[b]{$x=y$}}
\put(130,65){\makebox(0,0)[b]{$x=y$}}
\put(170,65){\makebox(0,0)[b]{$y=z$}}
\put(10,65){\makebox(0,0)[b]{$\top$}}
\put(150,65){\makebox(0,0)[b]{$\wedge$}}

\put(-5,40){\makebox(0,0)[r]{$r_x$}}
\put(55,40){\makebox(0,0)[r]{$s_{x,y}$}}
\put(125,40){\makebox(0,0)[r]{$t_{x,y,z}$}}

\end{picture}
\end{center}

Then in the category $\cal C$ of proofs defined via generality we find,
for example, the following equations:
\[
\begin{array}{l}
t_{x,y,y}\circ(1_{x=y}\wedge r_y)=\; \k{\wedge}{1}_{x=y,\top},
\\
t_{x,x,y}\circ(r_x\wedge 1_{x=y})=\; \k{\wedge}{2}_{\top,x=y},
\\
s_{y,x}\circ s_{x,y}=1_{x=y},
\\
s_{x,x}\circ r_x=r_x.
\end{array}
\]

\vspace{.2cm}

\section{Representing equivalence relations by sets of functions}

In this section we will consider some general matters
concerning the representation of arbitrary equivalence relations.
This will serve for demonstrating in the next section the isomorphism of
our representation of \G\ in the category whose arrows are
binary relations between finite ordinals.

Let $X$ be an arbitrary set, and let $R\subseteq X^2$.
Let $p$ be a set such that for $p_0 \neq p_1$ we have
$p_0,p_1\in p$. Let $S\subseteq p^2$, and for $i,j\in\{0,1\}$
let $p_i S p_j$ iff $i\leq j$. Consider the following set of functions:
\[
{\cal F}^S(R)=_{\mbox{\scriptsize\it def}} \{f:X\str p\mid
(\forall x,y\in X)(xRy\Rightarrow (f(x),f(y))\in S)\},
\]
and for $x,y\in X$ let $f_x:X\str p$ be defined as follows:
\[
f_x(y)=\left\{\begin{array}{ll}
p_1 & {\mbox{\rm if $xRy$}}
\\
p_0 & {\mbox{\rm if not $xRy$.}}
\end{array}
\right.
\]
The function $f_x$ is the characteristic function of the $R$-cone
over $x$. We can then prove the following proposition.

\vspace{.3cm}

\noindent {\sc Proposition} 1.\quad {\it The relation $R$ is transitive
iff $(\forall x\in X)f_x\in{\cal F}^S(R)$.}

\vspace{.2cm}

\noindent {\it Proof.}\quad $(\Rightarrow)$ Suppose $yRz$.
If $f_x(y)=p_0$, then $(f_x(y),f_x(z))\in S$. If $f_x(y)=p_1$, then
$f_x(z)=p_1$ by the transitivity of $R$.

\vspace{.1cm}

$(\Leftarrow)$ If $yRz\Rightarrow (f_x(y),f_x(z))\in S$, then
$yRz\Rightarrow (f_x(y)=p_0\;{\mbox{\it or}}\; f_x(z)=p_1)$, which
means $yRz\Rightarrow(xRy\Rightarrow xRz)$.
\qed

\vspace{.2cm}

With the definition
\[
{\cal F}^=(R)=_{\mbox{\scriptsize\it def}}\{f:X\str p\mid
(\forall x,y\in X)(xRy\Rightarrow f(x)=f(y))\}
\]
we can prove the following proposition.

\vspace{.3cm}

\noindent {\sc Proposition} 2.\quad {\it If $S$ is reflexive and
antisymmetric, and $R$ is symmetric, then ${\cal F}^=(R)={\cal F}^S(R)$.}

\vspace{.2cm}

\noindent {\it Proof.}\quad If $f\in {\cal F}^=(R)$,
then $f\in {\cal F}^S(R)$ by the
reflexivity of $S$; and if $f\in {\cal F}^S(R)$, then $f\in {\cal F}^=(R)$
by the symmetry of $R$ and the antisymmetry of $S$.
\qed

\vspace{.2cm}

Suppose now the relation $S$ is reflexive and antisymmetric.
We can then prove the following proposition.

\vspace{.3cm}

\noindent {\sc Proposition} 3.\quad {\it The relation $R$ is
an equivalence relation iff $(\forall x,y\in X)(xRy$ $\Leftrightarrow
(\forall f\in {\cal F}^=(R))f(x)=f(y))$.}

\vspace{.2cm}

\noindent {\it Proof.}\quad Note first that in
the equivalence on the right-hand
side the direction $xRy\Rightarrow (\forall f\in{\cal F}^=(R)) f(x)=f(y)$
is satisfied by definition.

\vspace{.1cm}

$(\Rightarrow)$
If $(\forall f\in {\cal F}^=(R))f(x)=f(y)$, then
$f_x(x)=f_x(y)$ by Propositions 1 and 2, and since
$f_x(x) = p_1$ by the reflexivity of $R$, we have $xRy$.

\vspace{.1cm}

$(\Leftarrow)$ From
\[(\ast)\quad (\forall x,y\in X)((\forall f\in {\cal F}^=(R)) f(x)=f(y)
\Rightarrow xRy)
\]
we obtain that $R$ is reflexive by taking that $x=y$.

For transitivity, suppose $xRy$ and $yRz$. Then for every
$f \in {\cal F}^=(R)$
we have $f(x)=f(y)=f(z)$, and hence $xRz$ by $(\ast)$.

For symmetry, suppose $xRy$. Then for every $f \in {\cal F}^=(R)$
we have $f(y)=f(x)$, and hence $yRx$ by $(\ast)$.
\qed

\vspace{.3cm}

\noindent {\sc Corollary.} \quad {\it If $R_1,R_2\subseteq X^2$ are
equivalence relations, then $R_1= R_2$ iff
${\cal F}^=(R_1)={\cal F}^= (R_2)$.}

\vspace{.2cm}

{\it Proof.}\quad We have, of course, that $R_1= R_2$ implies
${\cal F}^=(R_1)={\cal F}^= (R_2)$. The converse follows
from Proposition 3. \qed

\vspace{.2cm}

\section{Representing \G\ in the category \R}

Let \R\ be the category whose arrows are binary relations between finite
ordinals. Let $I_n \subseteq n\times n$ be the identity relation on $n$;
the composition $R_2\circ R_1\subseteq n\times k$
of $R_1\subseteq n\times m$ and $R_2\subseteq m\times k$ is
$\{(x,y)\mid (\exists z \in m)(xR_1z \;{\mbox {\it and}}\; zR_2y)\}$.

Then for $p\in \omega$ such that $p\geq 2$ we define a functor $F_p$
from \G\ to \R\ in the following manner. On objects $F_p$ is defined by
$F_p(n)=p^n$. Every element of the ordinal $p^n$ is identified by a
function $f:n\str p$ in the following way. Every $f:n\str p$
corresponds to a sequence of members of $p$ of length $n$. These
sequences can be ordered lexicographically, and $(p^n,\leq)$ is
isomorphic to the set of these sequences ordered lexicographically.

For $f_1: n\str p$ and $f_2:m\str p$, let $[f_1,f_2]:n\cup m^{+n}\str p$
be defined by
\[
[f_1,f_2](x)=\left\{\begin{array}{ll}
f_1(x) & {\mbox{\rm if }}x \in n
\\
f_2(x-n) & {\mbox{\rm if }}x \in m^{+n}.
\end{array}
\right.
\]

For $R:n\str m$ an arrow of \G, and for $f_1:n\str p$ and $f_2:m\str
p$, we define $F_p(R)$ by
\[(f_1,f_2)\in F_p(R)\quad{\mbox{\rm iff}}\quad [f_1,f_2]\in{\cal F}^=(R),
\]
where ${\cal F}^=(R)$ is the set of functions defined as in the preceding
section. Here $X$ is $n+m$, while for
the ordinal $p\geq 2$ we have that $p_0$ is $0$, $p_1$ is $1$ and
$S$ is $\leq$. It
remains to check that $F_p$ so defined is indeed a functor.

\vspace{.3cm}

\noindent {\sc Proposition} 4.\quad {\it $F_p$ is a functor.}

\vspace{.2cm}

\noindent {\it Proof.}\quad We show first that $F_p(\mj_n)=I_{p^n}$
by remarking that
\[
[f_1,f_2]\in {\cal F}^=(\mj_n)\quad{\mbox{\rm iff}}\quad f_1=f_2.
\]

Next we have to show that $F_p(R_2\ast R_1)=F_p(R_2)\circ F_p(R_1)$ for
$R_1:n\str m$ and $R_2: m\str k$. This amounts to showing that for
$f_1: n\str p$ and $f_2:k\str p$
\[
(\ast\ast)\quad (\forall x,y\in n\cup k^{+n})(x(R_2\ast R_1)y\Rightarrow
[f_1,f_2](x)=[f_1,f_2](y))
\]
is equivalent to the assertion that there exists an $f_3:m\str p$ such
that the following two statements are satisfied:
\[\begin{array}{l}
(\ast 1)\quad (\forall x,y\in n\cup m^{+n})(x R_1 y\Rightarrow
[f_1,f_3](x)=[f_1,f_3](y)),
\\[.1cm]
(\ast 2)\quad (\forall x,y\in m\cup k^{+m})(x R_2 y\Rightarrow
[f_3,f_2](x)=[f_3,f_2](y)).\end{array}
\]

First we show that $(\ast\ast)$ implies that for some $f_3:m\str p$ we have
$(\ast 1)$ and $(\ast 2)$. Before defining $f_3$ we check the
following facts that follow from $(\ast\ast)$:
\vspace{.2cm}
{\footnotesize
\[
\begin{array}{ll}
{\mbox {\rm (I)}} &
(\forall x,y\in n)((x,z),(y,z)\in\cl(R_1\cup R_2^{+n})\Rightarrow
f_1(x)=f_1(y)),
\\[.1cm]
{\mbox {\rm (II)}} &
(\forall x,y\in k^{+m+n})((x,z),(y,z)\in\cl(R_1\cup R_2^{+n})\Rightarrow
f_2(x-m-n)=f_2(y-m-n)),
\\[.1cm]
{\mbox {\rm (III)}} &
(\forall x\in n)(\forall y\in k^{+m+n})((x,z),(y,z)\in\cl(R_1\cup
R_2^{+n}) \Rightarrow f_1(x)=f_2(y-m-n)).
\end{array}
\]}

Then for $z\in m$ we define $f_3(z)\in p$ as follows. If for some
$x\in n\cup k^{+m+n}$ we have $(x,z+n)\in\cl(R_1\cup R_2^{+n})$, then
\[
f_3(z)=\left\{
\begin{array}{ll}
f_1(x) & {\mbox{\rm if }} x\in n
\\
f_2(x-m-n) & {\mbox{\rm if }} x\in k^{+m+n}.
\end{array}
\right.
\]
This definition is correct according to (I)-(III). Otherwise, if
there is no such $x$, we put $f_3(z)=0$.

Let us now demonstrate $(\ast 1)$. Suppose $xR_1 y$. If $x,y\in
n$, then $f_1(x)=f_1(y)$ by $(\ast\ast)$. If $x,y\in m^{+n}$, then
$f_3(x-n)=f_3(y-n)$ by the definition of $f_3$. If $x\in n$ and
$y\in m^{+n}$, then $f_1(x)=f_3(y-n)$ by the definition of $f_3$.
In every case we obtain $[f_1,f_3](x)=[f_1,f_3](y)$. We
demonstrate analogously $(\ast 2)$.

It remains to show that from the assumption that for some $f_3:m\str p$
we have $(\ast 1)$ and $(\ast 2)$ we can infer $(\ast\ast)$. Suppose
$x(R_2\ast R_1)y$ and $x,y\in n$. Then there is a sequence
$x_1,\ldots, x_{2l}$, $l\geq 1$, such that $x=x_1$, $y=x_{2l}$, for
every $i$ such that $1<i<2l$ we have $x_i\in m^{+n}$,
and
\[
x_1 R_1 x_2, \; x_2R_2^{+n} x_3, \; x_3 R_1 x_4,\ldots,
\; x_{2l-1} R_1 x_{2l}.
\]
Then by
applying $(\ast 1)$ and $(\ast 2)$ we obtain $f_1(x)=f_1(y)$,
and hence $[f_1,f_2](x)=[f_1,f_2](y)$. We proceed analogously when
$x,y\in k^{+n}$, or when $x\in n$ and $y\in k^{+n}$.
\qed

\vspace{.2cm}

Next we show that $F_p$ is faithful. Since $F_p$ is one-one on objects,
this amounts to showing that it is one-one on arrows.

\vspace{.3cm}

\noindent {\sc Proposition} 5.\quad {\it $F_p$ is faithful.}

\vspace{.2cm}

\noindent {\it Proof.}\quad Suppose $F_p(R_1)=F_p(R_2)$.
This means that for
every $f_1:n\str p$ and every $f_2: m\str p$ we have $[f_1,f_2]\in
{\cal F}^=(R_1)$ iff $[f_1,f_2]\in {\cal F}^=(R_2)$. But every function
$f:n+m\str p$ is of the form $[f_1,f_2]$ for some $f_1:n\str p$ and
some $f_2:m\str p$. Hence ${\cal F}^=(R_1)={\cal F}^=(R_2)$, and $R_1=R_2$
by the Corollary of the preceding section.
\qed

\vspace{.2cm}

So we have an isomorphic representation of \G\ in \R.

\section{Connection with Brauer algebras}

The representation of the category \G\ in the category \R,
which we have presented in the preceding section, is closely connected
to Brauer's representation of Brauer algebras, which is the orthogonal
group case of \cite{B37} (Section 5; see also \cite{W88}, Section 2,
\cite{J94}, Section 3, and \cite{DKP02}).

An $(n,n)$-{\it diagram} is an arrow $R:n\str n$ of \G\ such that
every member of the partition that corresponds to $R$ is a two-element
set. When defining the composition $R_2\ast R_1:n\str k$ of $R_1:n\str
m$ and $R_2:m\str k$ we didn't take into account the ``circles'', or
``closed loops'', in $\cl(R_1\cup R_2^{+n})$, namely, those members $X$
of the partition corresponding to $\cl(R_1\cup R_2^{+n})$ such that
$X\cap n=\emptyset$ and $X\cap k^{+m+n}=\emptyset$.
(In our example in Section 2 we have a circle involving 7 and 8
in the drawing.) Let $l(R_1,R_2)$ be
the number of those circles.

For $n\in \omega$ and $c$ a complex number,
the {\it Brauer algebra} $B(n,c)$ is the
algebra with basis the set of all $(n,n)$-diagrams, and multiplication
between two $(n,n)$-diagrams $R_1$ and $R_2$ defined as
$c^{l(R_1,R_2)}(R_2\ast R_1)$ (cf. \cite{J94}, Section 2).

Let $V$ be a vector space of dimension $p$ with basis $w_0,\ldots, w_{p-1}$.
For the $(n,n)$-diagram $R$ define $\beta(R)\in {\mbox{\rm
End}}(\otimes^n V)$ by the matrix (with respect to the basis
$\{w_{a_0}\otimes \ldots \otimes w_{a_{n-1}} \mid a_i\in\{0,\ldots,
p-1\}\}$ of $\otimes^n V$)
\[
\beta(R)_{a_0\ldots a_{n-1},a_n\ldots a_{2n-1}}
=_{\mbox{\scriptsize\it def}}
\prod_{\{i,j\}\in R}\delta(a_i,a_j)
\]
where $\delta$ is the Kronecker $\delta$ (cf. \cite{J94}, Definition
3.1). The sequences $a_0\ldots a_{n-1}$ and $a_n\ldots a_{2n-1}$
correspond to functions from $n$ to $p$. This defines a homomorphism of
the Brauer algebra $B(n,p)$ onto a subalgebra of ${\mbox{\rm
End}}(\otimes^n V)$ (see \cite{J94}, Lemma 3.2), which is Brauer's
representation of Brauer algebras in \cite{B37} (Section 5).

Every $n\times m$-matrix $M$ whose entries are only 0 and 1 may be
identified with a binary relation ${\cal R}_M\subseteq n\times m$ such that
$M(i,j)=1$ iff $(i,j)\in {\cal R}_M$. The matrix $M$ is the characteristic
function of ${\cal R}_M$. The $p^n\times p^n$-matrix $\beta(R)$ is a 0-1
matrix such that ${\cal R}_{\beta(R)}$ is precisely $F_p(R)$. So at the root of
Brauer's representation we find a particular case of our
representation of \G\ in \R.

One difference between our representation and Brauer's is that we deal
with arbitrary equivalence relations, whereas Brauer's equivalence
relations are more special---they correspond to partitions whose every
member is a two-element set. Another difference is that we have
composition of relations where Brauer has multiplication of matrices.
One passes from multiplication of matrices to our composition of
relations by disregarding ``circles'' and by taking that $1+1=1$.

\vspace{.5cm}

\noindent {\footnotesize {\it Acknowledgement.} The writing of this paper
was financed by the Ministry of Science, Technology and Development of
Serbia through grant 1630 (Representation of proofs with applications,
classification of structures and infinite combinatorics).}

\vspace{.5cm}

\noindent {\footnotesize {\it Acknowledgement added on 8 April 2016.} We would like to thank 
Sonja Telebakovi\' c for helping us to correct in this [v7] version of our paper a mistake concerning the function $f_3$ (in the definition of this function and in the demonstration of $(\ast 1)$ in the first paragraph after this definition) in Section~5 of the previous versions of the paper.}


\newpage

\begin{thebibliography}{99}

\bibitem{B37} R. Brauer, On algebras which are connected with the
semisimple continuous groups, {\it Ann. of Math.} 38 (1937), pp.
857-872.

\bibitem{D02} K. Do\v sen, Identity of proofs based on normalization and
generality, 2002 (available at: http:// arXiv. org/ math. LO/ 0208094).

\bibitem{DKP02} K. Do\v {s}en, \v Z. Kovijani\' c
and Z. Petri\'{c}, A new proof of the faithfulness of Brauer's
representation of Temperley-Lieb algebras, 2002 (available at:
http:// arXiv. org/ math. GT/ 0204214).

\bibitem{DP01} K. Do\v sen and Z. Petri\' c, The maximality of
cartesian categories, {\it Math. Logic Quart.} 47 (2001),
pp. 137-144 (available at: http:// arXiv. org/ math. CT/ 9911059).

\bibitem{DP02} K. Do\v {s}en and Z. Petri\'{c}, Bicartesian coherence,
{\it Studia Logica} 71 (2002), pp. 331-353 (available at:
http:// arXiv. org/ math. CT/ 0006052).

\bibitem{DP03} K. Do\v {s}en and Z. Petri\'{c}, A Brauerian
representation of split preorders, 2002 (available at:
http:// arXiv. org/ math. LO/ 00211277).

\bibitem{EK66} S. Eilenberg and G.M. Kelly,
A generalization of the functorial calculus, {\it J. Algebra}
3 (1966), pp. 366-375.

\bibitem{J94} V.F.R. Jones, A quotient of the affine Hecke algebra in the
Brauer algebra, {\it Enseign. Math. (2)} 40 (1994), pp. 313-344.

\bibitem{K95} C. Kassel, {\it Quantum Groups}, Springer, Berlin, 1995.

\bibitem{KML71} G.M. Kelly and S. Mac Lane, Coherence in closed categories,
{\it J. Pure Appl. Algebra} 1 (1971), pp. 97-140, 219.

\bibitem{L68} J. Lambek, Deductive systems and categories I: Syntactic
calculus and residuated categories, {\it Math. Systems Theory} 2 (1968),
pp. 287-318.

\bibitem{L69} J. Lambek, Deductive systems and categories II: Standard
constructions and closed categories, in: {\it Category Theory,
Homology Theory and their Applications I}, Lecture Notes in Math. 86,
Springer, Berlin, 1969, pp. 76-122.

\bibitem{L72} J. Lambek, Deductive systems and
categories III: Cartesian closed categories, intuitionist propositional
calculus, and combinatory logic, in: F.W. Lawvere ed., {\it Toposes,
Algebraic Geometry and Logic}, Lecture Notes in Math. 274, Springer,
Berlin, 1972, pp. 57-82.

\bibitem{L74} J. Lambek, Functional completeness of cartesian categories,
{\it Ann. Math. Logic} 6 (1974), pp. 259-292.

\bibitem{LS86} J. Lambek and P.J. Scott, {\it Introduction to Higher-Order
Categorical Logic}, Cambridge University Press, Cambridge, 1986.

\bibitem{P02} Z. Petri\' c, Coherence in substructural
categories, {\it Studia Logica} 70 (2002), pp. 271-296
(available at: http:// arXiv. org/ math. CT/ 0006061).

\bibitem{P71} D. Prawitz, Ideas and results in proof theory, in:
J.E. Fenstad ed., {\it Proceedings of the Second Scandinavian
Logic Symposium}, North-Holland, Amsterdam, 1971, pp. 235-307.

\bibitem{S75} M.E. Szabo, A counter-example to coherence in
cartesian closed categories, {\it Canad. Math. Bull.} 18 (1975),
pp. 111-114.

\bibitem{W88} H. Wenzl, On the structure of Brauer's centralizer
algebras, {\it Ann. of Math.} 128 (1988), pp. 173-193.

\end{thebibliography}
\end{document}